\input amstex
\magnification=1200
\font\cyr=wncyr10
\font\cyi=wncyi10
\font\cyb=wncyb10
\font\cyre=wncyr8
\font\cyie=wncyi8
\font\cybe=wncyb8
\documentstyle{amsppt}
\NoRunningHeads
\NoBlackBoxes
\define\soa{\operatorname{\frak s\frak o}}
\define\sla{\operatorname{\frak s\frak l}}
\define\sltwo{\sla(2,\Bbb C)}
\define\Cvir{\operatorname{\Bbb C\frak v\frak i\frak r}}
\define\crc{\Bbb S^1}
\define\Vect{\operatorname{Vect}(\crc)}
\define\CVect{\operatorname{\Bbb CVect}(\crc)}
\define\sothree{\soa(3)}
\define\sofour{\soa(4)}
\define\CA{\Cal A}
\define\CK{\Cal K}
\define\BC{\bold C}
\define\BO{\bold O}
\define\Fg{\frak g}
\define\Fp{\frak p}
\define\Fv{\frak v}
\define\Fw{\frak w}
\define\End{\operatorname{End}}
\define\Mor{\operatorname{Mor}}
\define\Abel{\operatorname{\Cal A\Cal B\Cal E\Cal L}}
\define\Add{\operatorname{\Cal A\Cal D\Cal D}}
\define\Lie{\operatorname{\Cal L\Cal I\Cal E}}
\define\LC{\operatorname{LC}}

\define\HS{\operatorname{\Cal H\Cal S}}

\topmatter
\title $q_R$--conformal symmetries in 2D nonlocal quantum field theory,
categorical representation theory, and Virasoro algebra
\endtitle
\author\bf D.V.Juriev\footnote{\ This is an English translation of the
original Russian version, which is located at the end of the article as an
appendix. In the case of any differences between English and Russian
versions caused by a translation the least has the priority as the original
one.\newline}
\endauthor
\date q-alg/9712009\enddate
\affil\eightpoint Research Center for Mathematical Physics and Informatics
``Thalassa Aitheria'',\linebreak
ul.Miklukho-Maklaya 20-180, Moscow 117437 Russia.\linebreak
E-mail: denis\@juriev.msk.ru
\endaffil
\abstract The paper is devoted to the symmetry aspects of 2D nonlocal
field theory, which is the simplest deformation of the conformally invariant
quantum field theory with one free bosonic field. The inverse
problem of representation theory is solved for $q_R$--conformal
symmetries, which are infinite dimensional hidden sym\-met\-ries of the
field theory. There is developed a mathematical formalism based on the
abstract categorical representation theory and sufficient for
applications of symmetry methods to the investigation of processes
of synthesis, decay and interaction of quasiparticles in the models
of 2D nonlocal quantum field theory.
\endabstract
\endtopmatter
\document
This paper being devoted to one of aspects of the inverse problem of
rep\-re\-sen\-tation theory [1] is addressed as to specialists in mathematical
physics and applied mathematics involved in the analysis of hidden
symmetries in quantum field models, classical and quantum dynamical and
controlled systems as to mathematicians-algebraists, who are interested
in new aspects of representation theory that are explicated during the
investigations of concrete models and systems. An exposition of the
material has its goals, in general, to show how unexpected mathematical
objects appear during the analysis of concrete hidden symmetries and
in what extend the modern rather abstract mathematical concepts
(category theory, Gro\-then\-dieck topology, theory of sheaves and topoi) may
be essential for the un\-der\-stan\-ding of symmetry aspects of problems of
mathematical physics. It is rather convenient to remind that an initiative
of the explication of such applications of modern
`categoric-algebrogeometric' apparatus to problems of theoretical and
mathematical physics belongs to Yu.I.Manin (see e.g.[2]).

The inverse problem of representation theory is to restore an abstract
ma\-the\-ma\-ti\-cal object (or its characteristics) by a concrete realization,
for example, by a fixed set of matrices or operators in an infinite
dimensional space [1]. Thus, one may correspond a Lie algebra
$\soa(n,\Bbb R)$ to all real $n\times n$ skew-symmetric matrices.
More complicated examples were considered, for example, in [1,3].

Often the inverse problem of representation theory means to search
a rep\-re\-sen\-ta\-tive from a fixed class of abstract algebraic structures
(e.g. Lie algebras or associative algebras with quadratic relations)
by its representation, it is defined at the same time what the
representation of an algebraic object of this class is. Under the
solution of such inverse problem as a rule new rather interesting
mathematical objects appear. In this way there were unraveled, for
example, quantum groups [4]; some other illustrations (Racah-Wigner
algebras, $\mho$--algebras, Sklyanin algebras etc.) are contained in
[1,3].

However, the inverse problem may be considered in another manner. Let us
suppose that an abstract object belongs to a rather simple class of
structures, for example, is a Lie algebra, but the form of representation
differs from a standard one. Such formulation of a problem was
proposed in [5]. Of course, methodologically the change of the
foreshortening means the change of the search of similarity between
the concrete and the abstract objects (the search of homomorphism
of the least into the first) to the analysis of a possibility of
a representation of the least by the first. That presupposes in any
case a question what a representation is and what `to represent' means.
Therefore, the using of the abstract algebraic apparatus of the category
theory is reasonable, that was done in [6]. Thus, this article,
discussing how the Virasoro algebra (more precisely, the Witt algebra,
whose central extension the Virasoro algebra is) appears under an
application of the categoric representation theory to analysis of
$q_R$--conformal symmetries in 2D nonlocal $q_R$--conformal field theory,
is a systematic exposition of a material of electronic preprints [5,6].

\head 1. $q_R$--conformal field theories and $q_R$--conformal
symmetries [7]
\endhead

The spectrum of fields of the simplest model of conformal field theory
is generated by one free bosonic field and, thus, consists of all
composite fields, which are the normal ordered pointwise
products of the current constructed from the free field (see a
general construction of currents from Fubini-Veneziano fields in
[7]) on itself. Other models of the conformal field theory, which
complete description was done in [8], are received under a transition
to several free fields (the so-called free field representation [9]).

The model with one free bosonic field is described by a Lagrangian
$$L=\int\varphi(z)\bar\varphi(\bar z)\,dzd\bar z,$$
where an integration is done inside the unit complex disk. From many points
of view (as theoretical as technical) it is useful to consider a nonlocal
deformation (cf.[10]) of the model with a nonlocal Lagrangian
$$L_h=\int K_h(z,\bar z)\varphi(z)\bar\varphi(\bar z)\,dzd\bar z,$$
where the kernel $K_h(z,\bar z)$ is the Bergman kernfunction of the
Lobachevski{\v\i} metric in the unit complex disk:
$$K_h(z,\bar z)=(1-z\bar z)^{-2-2h},$$
so the nonlocal Lagrangian is invariant under all M\"obius
(linear-fractional) trans\-for\-mations of the complex disk. Note that
the Lagrangian of the free bosonic field $L(\varphi)$ may be obtained
from the nonlocal Lagrangians $L_h$ as their limit at $h\to\infty$.
Models of 2D quantum field theory with the Lagrangian $L_h$ (and whose
classical counterparts have the natural nonlocal nondegenerate Poisson
brackets with Lagrangian $L_h$ as a K\"ahler
potential of the K\"ahler metric constructed from them) are the simplest
models of nonlocal theory and at the same time they possess an interesting
and nontrivial internal algebraic structure.

A language of the operator formalism of the quantum field theory
was applied to these nonlocal models after their
holomorphic-antiholomorphic (chiral) fac\-to\-ri\-zation in the author's
paper [7], where the models were called $q_R$--conformal field theories
($q_R=\frac1{2h-1}$). The components of operator fields of the
$q_R$--conformal field theory, which are invariant under all complex
projective transformation of the Riemann sphere, admit a realization
by tensor operators in the Verma modules over the Lie algebra
$\sltwo$ and, therefore, may be written as (generalized) differential
operators, so many problems of the field theory for these models have
analytic or numeric solutions and the formulation of their difference
(lattice) approximations does not produce any problems (for example,
the problem of integration of a quantum-field top, which was considered
for the conformal field theory in [11] and for $q_R$--conformal theories
in [12]).

An important feature of 2D conformal field theories is the presence
of infinite additional symmetries, which are described by the Witt
algebra or by its central extension, the Virasoro algebra $\Cvir$ [13].
These symmetries appear as components of the decomposition of the
(chirally factorized) stress-energy tensor, the square of the free
bosonic field (more precisely, of the current constructed from it).
Analogs of the free field representations and stress-energy tensor
in the $q_R$--conformal case were considered in detail in the article [7].
Thus, the $q_R$--conformal stress-energy tensor is the generating function
of the spin 2 tensor operators in the Verma module over the Lie algebra
$\sltwo$. If the Verma module $V_h$ (here $h$ is the extremal weight [14])
over this algebra is realized in the space of all polynomials $\Bbb C[z]$
of a complex variable $z$, and the action of generators $L_i$ ($i=-1,0,1$;
$[L_i,L_j]=(i-j)L_{i+j}$) of the Lie algebra $\sltwo$ in it is written as
$$L_{-1}=z,\quad L_0=z\partial_z+h,\quad L_1=z\partial_z^2+2h\partial_z,$$
then the components of the $q_R$--conformal stress-energy tensor
($q_R$--conformal sym\-met\-ries) are of the form:
$$ L_{-k}=(\xi+h(k+1))\partial_z^k\quad (k\ge0),\quad
L_k=z^k\frac{\xi+h(k+1)}{(\xi+2h)\ldots(\xi+2h+k-1)}\quad (k\ge1),$$
where $\xi=z\partial_z.$

The $q_R$--conformal symmetries form an infinite family. There
were done several attempts to solve the inverse problem of representation
theory for them, i.e. to unravel their algebraic structure. One of these
attempts was described in the article [7], its result is the
associative algebra $\Cvir(q_R)$, an analog of the nonlinear $\sla_2$ of
[15]. Another approach was sketched in the electronic preprint [5]
and developed in the author's article `Approximate representations
and the Virasoro algebra' prepared for a publication. However, there is
a necessity of the further search of other ways to solve the inverse
problem of representation theory for $q_R$--conformal symmetries.
The next paragraph id devoted to a description of one of the possible
versions.

\head 2. Algebraic structure of $q_R$--conformal symmetries
\endhead

\definition{Definition 1 [5]}

{\bf A.} A linear space $\Fv$ is called {\it a Lie composite\/}
iff there are fixed its subspaces $\Fv_1,\ldots \Fv_n$ ($\dim\Fv_i>1$)
supplied by the compatible structures of Lie algebras. Com\-pa\-ti\-bi\-li\-ty means that
the structures of the Lie algebras induced in $\Fv_i\cap\Fv_j$ from $\Fv_i$ and
$\Fv_j$ are the same. The Lie composite is called {\it dense\/} iff
$\Fv_1\uplus\ldots\uplus\Fv_n=\Fv$ (here $\uplus$ denotes the sum of linear
spaces). The Lie composite is called {\it connected\/}
iff for all $i$ and $j$ there exists a sequence $k_1,\ldots k_m$ ($k_1=i$,
$k_m=j$) such that $\Fv_{k_l}\cap\Fv_{k_{l+1}}\ne\varnothing$.

{\bf B.} {\it A representation\/} of the Lie composite $\Fv$ in the space $H$
is the linear mapping $T:\Fv\mapsto\End(H)$ such that $\left.T\right|_{\Fv_i}$
is a representation of the Lie algebra $\Fv_i$ for all $i$.

{\bf C.} Let $\Fg$ be a Lie algebra. A linear mapping $T:\Fg\mapsto\End(H)$
is called {\it the composed representation\/} of $\Fg$ in the linear space
$H$ iff there exists a set $\Fg_1,\ldots,\Fg_n$ of the Lie subalgebras of
$\Fg$, which form a dense connected composite and $T$ is its representation.
\enddefinition

Reducibility and irreducibility of representations of the Lie composites are
de\-fined in the same manner as for Lie algebras. One may also formulate a
superanalog of the Definition 1. The set of representations of the fixed Lie
composite is closed under the tensor product and, therefore, may be supplied by
the structure of {\it a tensor category\/} [16].

\example{Example 1 (The Octahedron Lie Composite) [5]} Let us consider an
oc\-ta\-hed\-ron with the vertices
$A$, $B$, $C$, $D$, $E$, $F$, the edges $(AB)$, $(AC)$, $(AD)$, $(AE)$,
$(BC)$, $(BF)$, $(CD)$, $(CF)$, $(DE)$, $(DF)$, $(EF)$, and the faces
$(ABC)$, $(ACD)$, $(ADE)$, $(AEB)$, $(BCF)$, $(CDF)$, $(DEF)$, $(EBF)$.
Let $\Fv$ be a six--dimensional linear space with the basis labelled by the
vertices of the octahedron, $\Fv_1$, $\Fv_2$, $\Fv_3$, $\Fv_4$ be four
three--dimensional subspaces in $\Fv$ corresponded to the faces
$(ABC)$, $(ADE)$, $(CDF)$, $(EBF)$. All subspaces $\Fv_i$ are supplied by the
structures of the Lie algebras isomorphic to $\sothree$ (such structures are
compatible to the orientations on the faces). The pentuple
$(\Fv,\Fv_1,\Fv_2,\Fv_3,\Fv_4)$ is a dense connected Lie composite.

\proclaim{Proposition} Let $T$ be an arbitrary representation of the Lie
composite $(\Fv;\Fv_1,\Fv_2,\mathbreak\Fv_3,\Fv_4)$ in the finite--dimensional
linear space $H$, then $H$ admits a representation of the Lie algebra $\sofour$.
If $T$ is an irreducible representation then there exist the real numbers
$\lambda_A$, $\lambda_B$, $\lambda_C$, $\lambda_D$, $\lambda_E$, $\lambda_F$
such that the operators $T(A)-\lambda_A\boldkey 1$, $T(B)-\lambda_B\boldkey 1$,
$T(C)-\lambda_C\boldkey 1$, $T(D)-\lambda_D\boldkey 1$, $T(E)-\lambda_E\boldkey
1$, $T(F)-\lambda_F\boldkey 1$ form an irreducible representation of $\sofour$.
\endproclaim

\demo{Proof} First, note that the commutator of operators corresponded to
the opposite vertices commute with operators corresponded to other four
vertices. It commutes with all six operators because they may be expressed as
commutators of the least four operators. So the commutator of operators
corresponded to the opposite vertices belongs to the center of the Lie
algebra generated by the all six operators. Let us factorize this Lie algebra
$\Fg$ by the center. Such quotient is isomorphic to $\sofour$ (one uses the
fact that formulas for commutators of all six operators are known up to the
center of $\Fg$). The statement of the theorem is a consequence of this result
and the fact that any central extension of the semisimple Lie algebra is
trivial (i.e. may be splitted -- see f.e.[17])
\qed\enddemo
\endexample

The construction of the Octahedron Lie composite may be generalized on
the certain class of polyhedra. However, any analogs of the Proposition
are not known for such general case.

\example{Example 2 (The Witt composite) [5]} Let $\Fw$ be the so--called Witt
algebra, which is a subalgebra of the complexification $\CVect$ of the Lie
algebra $\Vect$ of the smooth vector fields on a circle $\crc$ [18]. The Witt
algebra $\Fw$ consists of all polynomial vector fields and admits a basis
$e_k$ ($k\in\Bbb Z$) with commutation relations $[e_i,e_j]=(i-j)e_{i+j}$.
The Virasoro algebra $\Cvir$ [19] is one-dimensional nontrivial central
extension of the Witt algebra.

Let us consider two subalgebras $\Fp_{\pm}$ of $\Fw$ generated by $e_i$ with
$i\ge-1$ and $i\le1$; note that $\Fp_+\cap\Fp_-=\sltwo$. The triple
$(\Fw;\Fp_+,\Fp_-)$ is a dense connected Lie composite.

Each representation of the Virasoro algebra defines a representation of
the Witt composite in the same space. The reciprocal statement is not
correct, of course, what is shown by the following theorem.

\proclaim{Theorem 1A} The action of the Lie algebra $\sltwo$ in any Verma
module $V_h$ ($h$ is the highest weight) may be extended to the representation
of the Lie composite $(\Fw;\Fp_+,\Fp_-)$ and, hence, to the composed
representation of the Witt algebra $\Fw$.
\endproclaim

\demo{Proof} The Theorem follows from the explicit formulas for the spin 2
tensor ope\-ra\-tors in the Verma modules over the Lie algebra $\sltwo$
($q_R$--conformal symmetries), which realize a representation of the
Witt composite.\qed\enddemo
\endexample

\remark{Remark 1} The construction of the Witt composite may be generalized
on the Riemann surfaces of higher genus in lines of I.M.Krichever and
S.P.Novikov [19].
\endremark

\remark{Remark 2} Generalizing the terminology of [1,3] one may say that
the spin 2 tensor operators (i.e. the $q_R$--conformal symmetries)
in the Verma modules $V_h$ over the Lie algebra $\sltwo$
form the set of hidden symmetries, whose algebraic structure is one of the
Witt composite.
\endremark

Note that if the hidden symmetries realize a representation of the
Lie composite they should not be unpacked (a similar situation appears also
in the case of the isocommutator algebras of hidden symmetries and the
related Lie $\Fg$--bunches [3:To\-pic 3;1:\S2.1]).

\example{Example 3} Let $\Fw$ be the Witt algebra and $(\Fw;\Fp_{\pm})$ be the
Witt composite. Let us consider the abelian extension $\Fw^e$ of the Witt
algebra by the generators $f_i$ ($i\in\Bbb Z$) such that $[e_i,f_j]=jf_j$.
The subalgebras $\Fp_{\pm}$ of $\Fw$ may be extended to the subalgebras
$\Fp^e_{\pm}$ of $\Fw^e$ by the generators $f_i$, where $i\ge 0$ and $i\le 0$,
respectively. The triple $(\Fw^e;\Fp^e_{\pm})$ form the extended Witt composite.

\proclaim{Theorem 1B} The representation of the Witt composite in any Verma
module $V_h$ ($h$ is the highest weight) over the Lie algebra $\sltwo$ may be
extended to the representation of the Lie composite $(\Fw^e;\Fp^e_{\pm})$ and,
hence, the composed rep\-re\-sen\-tation of $\Fw$ in $V_h$ may be extended to the
composed rep\-re\-sen\-tation of $\Fw^e$.
\endproclaim

\demo{Proof} The additional generators $f_i$ are represented by the tensor
operators of spin 1 (the components of a decomposition of the $q_R$--affine
current [7]), namely, $f_i\mapsto\partial_z^i$ ($i\ge 0$),
$f_{-i}\mapsto z^i\frac1{(\xi+2h)\ldots(\xi+2h+i-1)}$ ($i\ge 1$)
\qed\enddemo
\endexample

\remark{Remark 3} It is very interesting to consider the composed
representations of the real semisimple Lie algebras $\Fg$, which unduce
representations of some natural sub\-al\-geb\-ras (for instance, of two opposite
maximal parabolic sub\-al\-geb\-ras or two opposite Borel sub\-al\-geb\-ras, perhaps plus
some $\sltwo$ imbed into $\Fg$, etc.).
\endremark

\remark{Remark 4} Composed representations of other algebraic structures
(isotopic pairs) were considered in the author's article [20].
\endremark

Thus, we established in the second paragraph that the $q_R$--conformal
symmetries may be regarded as a realization of the same algebraic
structure (the Witt algebra or the Virasoro algebra) as for the
conformal symmetries, however, with a new understanding of the
operation of `representation' itself. Thus, it is natural to formulate
a question what `to represent' means and what feautures should
the nonstandard representation have to allow to operate with them as
with standard ones, and how to receive new unknown types of
representations from the constructed ones. All these questions are
a subject of the following paragraph, which will have a rather formal
mathematical character being far from the initial formulation of the problem
from the nonlocal quantum field theory.

\head 3. Categorical representation theory and $q_R$--conformal
symmetries: composite and overlay representations
\endhead

\subhead 3.1. Elements of the categorical representation theory [6]
\endsubhead
We shall con\-si\-der the representations of classes of objects, which
constitute a category, which will be called {\it the ground category}.
The categorical aspects of the standard rep\-re\-sen\-tation theory were discussed
in [16]. Some categorical generalizations were described in [21]. However, we
shall formulate the most abstract settings, which are necessary for our
purposes.

\definition{Definition 2A} {\it A representation theory\/} for the
ground category $\CA$ is a con\-t\-ra\-va\-riant functor $R$ from the category
$\CA$ to the category $\Abel$ of all small abelian categories.
\enddefinition

Sometimes one should consider the category $\Add$ of all small additive
categories instead of $\Abel$. However, we shall consider the least
category for simplicity.

Often the ground category has some good properties, e.g. that for any
finite family of objects there exists their coproduct, which coincides
with their product. Such situation is realized for Lie algebras, Lie groups,
finite groups, associative algebras, Hopf algebras and many other structures.
However, the isotopic pairs (see f.e.[1:\S2.2;19]) and the most of other
algebraic pairs do not form a category of such type. For the ground
category $\CA$, in which products and coproducts of finite number of objects
exist and coincide, we shall claim in the definition of representation theories
that an associative family of imbeddings $R(a)\times R(b)\hookrightarrow
R(a+b)$ ($a$ and $b$ are any objects of the ground category $\CA$) is defined.
Such representation theories will be called {\it quasitensorial}.

\remark{Remark 5} If an object $a$ of the ground category $\CA$ admits a
coassociative mo\-no\-mor\-phism $\varepsilon$ into $a+a$ then $R(a)$ is
a tensor category iff the representation theory $R$ is quasitensorial.
\endremark

There exist non--quasitensorial representation theories even for the
well--known categories of the represented objects, e.g. general
$\HS$--projective representations of Lie algebras [5] or unitary
$\HS$--pseudorepresentations of Lie groups are out of this class.

\definition{Definition 2B} A representation theory for the ground category
$\CA$ is called {\it ho\-mo\-mor\-phic\/} iff there exists a subcategory $\CA_0$ of
$\CA$ ({\it the target subcategory\/}) such that for any object $a$ of $\CA$
the category $R(a)$ may be identified with the category $\Mor(a,\CA_0)$ of all
(equivalence classes of) morphisms from $a$ to the objects of the category
$\CA_0$.
\enddefinition

For instance, theories of all linear, projective, unitary representations of
Lie groups are homomorphic. Note that the target category $\CA_0$ is always
an additive subcategory of the ground category $\CA$.

\definition{Definition 2C} A representation theory for the ground
category $\CA$ is called {\it hid\-den\-ly homomorphic\/} iff there exists a
homomorphic representation theory $R'$ for a category $\CK$ and a functor
(multi-valued in general)
$\varrho:\CA\mapsto\CK$ such that $R=R'\circ\varrho$.
\enddefinition

Below we shall consider some examples and general constructions of the
hiddenly homomorphic representation theories for the ground category $\Lie$ of
the Lie al\-geb\-ras inspired by the considered above solution of the inverse
proble, of representation theory for $q_R$--conformal symmetries in
nonlocal quantum field theory, which are not homomorphic, and describe
their interpretations in terms of the categorical representation theory.

\subhead 3.2. Composed representation theories [6]
\endsubhead
Let us formulate an abstract categorical setting for the construction
of the composed representations defined in the second paragraph.

\definition{Definition 3A} Let $\CA$ be a topologized ground category
(i.e. supplied by a struc\-ture of the Grothendieck topology [22]). Let $R$
be a representation theory for $\CA$. {\it The composed representation
theory\/} $\BC(R)$ for $\CA$ may be constructed in the following manner.
Let $a$ be an object of the ground category $\CA$ and $S=(s_1,s_2,\ldots s_n)$
($s_i\in\Mor(a_i,a)$) be a cover of $a$ then the objects of the category
$\BC(R)(a)$ consists of all data $(b_1,b_2,\ldots b_n)$, $b_i\in R(a_i)$
such that for any object $c$ and monomorphisms $f\in\Mor(c,a)$ and
$f_i\in\Mor(c,a_i)$ ($f=s_i\circ f_i$) the equality ({\it the composite
glueing rule\/})
$$R(f_i)^*(b_i)=R(f_j)^*(b_j)$$
holds. The morphisms in $\BC(R)(a)$ are defined in the same manner.
\enddefinition

For any representation theory $R$ the composite representation theory
$\BC(R)$ is a sheaf of abelian categories over the topologized ground category
$\CA$ [22]. It is a sheaf canonically constructed from the pre--sheaf $R$
over the topologized ground category $\CA$ (note that the representation theory
for the topologized ground category $\CA$ is just a pre--sheaf over it).

\proclaim{Theorem 2A [6]} The composed representations of Lie algebras form a
composed representation theory $\BC(R)$, where $R$ is a standard
representation theory of Lie algebras (the covers of the Lie algebras are
defined by the dense connected Lie composites).
\endproclaim

Note that the Grothendieck topology of the Theorem 2A differs from the usual
one [21].

\remark{Remark 6} If $R$ is the standard representation theory then the
theory $\BC(R)$ is hiddenly homomorphic, the category $\CK$ is one of the
Lie composites, the category $\CK_0$ consists of Lie algebras $\End(H)$ for
all linear spaces $H$. i.e. just the same as for a homomorphic standard
represntation theory. However, if $R$ is a general representation theory
$\BC(R)$ is not obligatory hiddenly homomorphic.
\endremark

I suspect that the concept of the hidden homomorphicity of the composite
representation theories may be somehow understood in terms of the topos
theory [22].

\remark{Remark 7} $\BC(\BC(R))=\BC(R)$.
\endremark

\subhead 3.3. Overlay representation theories [6]
\endsubhead The disadvantages of the composed representation theory is clearly
explicated on the examples of the composed rep\-re\-sen\-tations of the Witt algebras
by the hidden infinite dimensional ($q_R$--conformal) symmetries in the Verma
modules over the Lie algebra $\sltwo$. First, the tensor product of a finite
number of these irreducible composed representations is ir\-re\-du\-cible.
This fact contradicts to the na{\"\i}ve intuition. Second, the hidden symmetries
do not form any representation themselves whereas intuitively they should form
the adjoint representation. From the point of view of the mathematical
physics all this disadvantages are essential because they do not allow to
use effectively the apparatus of representation theory for the analysis
of processes of decay, synthesis and interaction of quasiparticles in 2D
nonlocal field theory.

So one needs some generalization of the composed representations.
Let us define the operator Lie composites $\LC(H)$ as the sets of subspaces
$\End(H_i)$ in the spaces $\End(H)$ ($H=H_1+\ldots+H_m$) with the natural
structures of Lie algebras.

\definition{Definition 4}

{\bf A.} {\it An overlay representation\/} of the Lie composite $\Fv$ in the
space $H$ is the homomorphism $T$ of $\Fv$ into the operator Lie composite
$\LC(H)$.

{\bf B.} Let $\Fg$ be a Lie algebra. A linear mapping $T:\Fg\mapsto\End(H)$
is called {\it the overlay composed representation\/} (or simply {\it overlay
representation\/}) of $\Fg$ in the linear space $H$ iff there exists a set
$\Fg_1,\ldots,\Fg_n$ of the Lie subalgebras of $\Fg$, which form a dense
connected composite and $T$ is its overlay representation.
\enddefinition

\remark{Remark 8} The overlay representations of any Lie algebra $\Fg$
form a tensor category.
\endremark

The overlay representations solve the previously described difficulties.

\proclaim{Theorem 1C [6]} The tensor operators of spin 2 in the Verma modules
$V_h$ over the Lie algebra $\sltwo$ form an overlay representation of the Witt
algebra, which are subrepresentations of $\End(V_h)$.
\endproclaim

\remark{Remark 9} The tensor operators of any natural spin $n$ in the Verma
modules $V_h$ over the Lie algebra $\sltwo$ (described completely in [7]) form
overlay representations of
the Witt algebra, which are subrepresentations of $\End(V_h)$.
\endremark

Let us formulate the natural categorical settings for the construction
of overlay representations.

\definition{Definition 3B} Let $\CA$ be a topologized ground category. Let $R$
be a homomorphic representation theory for $\CA$ with the target subcategory
$\CA_0$ supplied by the Gro\-then\-dieck topology induced from $\CA$. {\it
The overlay representation theory\/} $\BO(R)$ for $\CA$ may be constructed
in the following manner. Let $a$ be an object of the ground category $\CA$ and
$S=(s_1,s_2,\ldots s_n)$ ($s_i\in\Mor(a_i,a)$) be a cover of $a$ then the
objects of the category $\BO(R)(a)$ consists of all data $(r_1,r_2,\ldots
r_n)$, $r_i\in\Mor(a,b_i)$, $b_i$ are objects of the target subcategory
$\CA_0$, which form a cover of the object $b$ of the same subcategory by the
monomorphisms $t_i\in\Mor(b_i,b)$, such that for any subobject $(c;p]$ of $a$
($p\in\Mor(c,a)$) the equality ({\it the overlay glueing rule\/})
$$r_i((a_i;s_i]\cap(c;p])\cap(b_j;t_j]=(b_i;t_i]\cap r_j((a_j;s_j]\cap(c;p])$$
holds. The morphisms in $\BO(R)(a)$ are defined in the same manner.
\enddefinition

However, I do not know a definition of the overlay representation theory
$\BO(R)$ for the representation theory $R$, which is not homomorphic.
Note that $\BO(R)$ is not a sheaf of abelian categories over $\CA$ in
general, and I do not know an abstract sheaf theoretical characterization
of the overlay representation theories.

\remark{Remark 10} The overlay representation theories $\BO(R)$ being defined
for the ho\-mo\-mor\-phic representations theories $R$ are hiddenly homomorphic.
\endremark

\proclaim{Theorem 3B [6]} The overlay representations of Lie algebras form
an overlay representation theory $\BO(R)$, where $R$ is a standard
representation theory of Lie algebras (the covers of the Lie algebras are
defined by the dense connected Lie composites and the target subcategory
$\CA_0$ consists of all Lie algebras $\End(H)$).
\endproclaim

Note that the Grothendieck topology of the Theorem 3B differs from the usual
one.

\remark{Remark 11} (A)
If $R$ is a homomorphic representation theory for the ground ca\-te\-go\-ry
$\CA$ then for any object $a$ of $\CA$ the category $\BC(R)(a)$ is a
subcategory of $\BO(R)(a)$.
(B) $\BC(\BO(R))=\BO(R)$.
\endremark

Thus, the symmetry aspects of 2D nonlocal field theory, which is
the simplest deformation of conformally invariant theory with one
free bosonic field, are con\-si\-dered in the article (\S1).
The inverse problem of representation theory is solved for
$q_R$--conformal symmetries, which are the infinite dimensional
hidden symmetries of the theory (\S2). Based on the abstract categorical
representation theory the mathematical apparatus sufficient for
applications of symmetry methods to in\-ves\-ti\-gations of processes of
synthesis, decay and interaction of quasiparticles in the models of
2D nonlocal quantum field theory is developed (\S3).

\Refs
\roster
\item" [1]" {\it Juriev D.V.}, An excursus into the inverse problem of
representation theory [in Russian]: mp\_arc/96-477.
\item" [2]" {\it Manin Yu.I.}, Gauge fields and complex geometry [in
Russian]. Moscow, Nauka, 1986.
\item" [3]" {\it Juriev D.}, Topics in hidden symmetries. I-IV,VI:
hep-th/9405050, q-alg/9610026, q-alg/9611003, q-alg/9611019, q-alg/9708028;
On the infinite-dimensional hidden sym\-met\-ries. I-III: funct-an/9612004,
funct-an/9701009, funct-an/9702002.
\item" [4]" {\it Reshetikhin N.Yu., Tahtadzhyan L.A., Faddeev L.D.},
St.Petersburg Math.J. 1 (1990) 193-225; {\it Drinfeld V.G.},  Zapiski
Nauchn.Semin.LOMI [in Russian] 155 (1986) 19-49; {\it Isaev A.P.},
Elem.Part.Atom.Nuclei [in Russian] 26(5) (1995) 1204-1263.
\item" [5]" {\it Juriev D.}, Topics in hidden symmetries. V:
funct-an/9611003.
\item" [6]" {\it Juriev D.}, Hidden symmetries and categorical
representation theory: q-alg/9612026.
\item" [7]" {\it Juriev D.V.}, Theor.Math.Phys. 101 (1994) 1387-1403.
\item" [8]" {\it Juriev D.V.}, St.Petersburg Math.J. 3 (1992) 679-686.
\item" [9]" {\it Morozov A.Yu.}, Elem.Part.Atom.Nuclei [in Russian]
23(1) (1992) 174-238.
\item"[10]" Nonlocal, nonlinear and nonrenormalizable field theories [in
Russian]. Dubna, JINR, 1970.
\item"[11]" {\it Alekseev A., Faddeev L., Semenov-Tian-Shansky M.},
Commun.Math.Phys. 149 (1992) 335-345.
\item"[12]" {\it Juriev D.V.}, Theor.Math.Phys. 106 (1996) 276-290.
\item"[13]" {\it Belavin A.A., Polyakov A.M., Zamolodchikov A.B.},
Nucl.Phys. B241 (1984) 333-380.
\item"[14]" {\it Zhelobenko D.P.}, Representations of the reductive Lie
algebras [in Russian]. Moscow, Na\-u\-ka, 1994.
\item"[15]" {\it Ro\v cek M.}, Phys.Lett. B255 (1991) 554-557.
\item"[16]" {\it Kirillov A.A.}, Elements of representation theory.
Springer, 1976; An introduction into the representation theory and
noncommutative harmonic analysis [in Russian]. Current Math.Problems.
Basic Directions 22, Moscow, VINITI, 1988.
\item"[17]" {\it Guichardet A.}, Cohomologie dea groupes topologiques et
des alg\`ebres de Lie. Cedic/Fer\-nan Nathan, Paris, 1980.
\item"[18]" {\it Fuchs D.B.}, Cohomology of infinite dimensional Lie
algebras [in Russian]. Moscow, Nau\-ka, 1984.
\item"[19]" {\it Krichever I.M., Novikov S.P.}, Funct.Anal.Appl. 21
(1987) 126-141; 21 (1987) 294-307; 23 (1989) 19-32; J.Geom.Phys. 5 (1988)
631-661 [reprinted in ``Geometry and physics. Essays in honour of
I.M.Gelfand'', Eds. S.Gindikin and I.M.Singer, Pitagora Editrice,
Bologna and Elsevier Sci.Publ., Amsterdam, 1991].
\item"[20]" {\it Juriev D.V.}, Theor.Math.Phys. 105 (1995) 1201-1209;
111 (1997) 511-518.
\item"[21]" {\it Geronimus A.Yu.}, Funkts.anal.i ego prilozh [in Russian].
5(3) (1971) 22-31.
\item"[22]" Th\'eorie des topos et cohomologie \'etale des sch\'emas.
Berlin, Heidelberg, New York, 1972-1973; Cohomologie \'etale. Berlin,
Heidelberg, New York, 1977.
\endroster
\endRefs

\ \newpage

\centerline{\bf APPENDIX: THE ORIGINAL RUSSIAN VERSION OF ARTICLE}
\ \newline
\centerline{\cyb $q_R$--KONFORMNYE SIMMETRII V DVUMERNOI0}
\centerline{\cyb NELOKALP1NOI0 KVANTOVOI0 TEORII POLYa, KATEGORNAYa}
\centerline{\cyb TEORIYa PREDSTAVLENII0 I ALGEBRA VIRASORO}
\ \newline
\centerline{\cyb D.V.Yurp1ev}
\ \newline\eightpoint
\centerline{\cyre Tsentr matematicheskoe0 fiziki i informatiki ``Talassa
E1teriya'',}
\centerline{\cyre ul.Miklukho-Maklaya 20-180, Moskva 117437 Rossiya.}
\centerline{E-mail: denis\@juriev.msk.ru}
\ \tenpoint\newline
\centerline{q-alg/9712009}
\ \newline

\eightpoint\cyre\ \newline
{\cybe Rezyume. }
V rabote rssmotreny simmetrie0nye aspekty dvumernoe0 nelokalp1noe0
te\-o\-rii polya, yavlyayushchee0sya prostee0shee0 deformatsiee0
konformno-invariantnoe0 kvan\-to\-voe0 teorii polya s odnim svobodnym
bozonnym polem. Reshena obratnaya zadacha teorii pred\-stav\-le\-nie0
dlya $q_R$--konformnyh simmetrie0, yavlyayushchih\-sya beskonechnomernymi
skry\-ty\-mi simmetriyami ukazannoe0 teorii. Na baze abstraktnoe0
kategornoe0 teorii pred\-s\-tav\-le\-nie0 razvit matematicheskie0 apparat,
dostatochnye0 dlya primeneniya sim\-met\-rie0\-nyh metodov k izucheniyu
protsessov sinteza, raspada i vzaimodee0stviya kvazichastits v
modelyah dvumernoe0 nelokalp1noe0 kvantovoe0 teorii polya.

\tenpoint\ \newline

\ \newline

\cyr
Dannaya rabota, posvyashchennaya odnomu iz aspektov obratnoe0
zadachi te\-o\-rii predstavlenie0 [1], adresovana kak spetsialistam v
matematicheskoe0 fi\-zi\-ke i prikladnoe0 matematike, zanimayushchihsya
analizom skrytyh sim\-met\-rie0 v kvantovo-polevyh modelyah, klassicheskih
i kvantovyh di\-na\-mi\-ches\-kih i upravlyaemyh sistemah, tak i
matematikam-algebraistam, in\-te\-re\-su\-yu\-shchim\-sya tem, kakie
novye storony teorii predstavlenie0 vskryvayut\-sya pri issledovanii
konkretnyh modelee0 i sistem. Izlozhenie materiala presleduet, v osnovnom,
tselp1 pokazatp1, kak pri analize konkretnyh skry\-tyh simmetrie0
voznikayut neozhidannye matematicheskie obp2ekty i v ka\-koe0 stepeni
sovremennye dostatochno abstraktnye matematicheskie kon\-tsep\-tsii
(teoriya kategorie0, topologii Grotendika, teoriya puchkov i toposov)
mogut bytp1 sushchestvenny pri osmyslenii simmetrie0nyh aspektov zadach
ma\-te\-ma\-ti\-ches\-koe0 fiziki. V dannoe0 svyazi vpolne umestno
napomnitp1, chto ini\-tsi\-a\-ti\-va raskrytiya podobnyh prilozhenie0
sovremennogo ``ka\-te\-gor\-no--al\-geb\-ro\-geo\-met\-ri\-ches\-ko\-go''
apparata k voprosam teoreticheskoe0 i ma\-te\-ma\-ti\-ches\-koe0 fiziki
prinadlezhit Yu.I.Maninu (sm.napr.[2]).

Obratnaya zadacha teorii predstavlenie0 zaklyuchaet\-sya v
vosstanovlenii abstraktnogo matematicheskogo obp2ekta (ili ego
harakteristik) po kon\-k\-ret\-noe0 realizatsii, naprimer, po zadannoe0
sovokupnosti matrits ili operatorov v beskonechnomernom prostranstve
[1]. Tak sovokupnosti vseh veshchestvennyh kososimmetricheskih
matrits $n\times n$ mozhno sopostavitp1 al\-geb\-ru Li $\soa(n,\Bbb R)$.
Bolee slozhnye primery, rassmatrivalisp1, naprimer, v [1,3].

Chasto obratnaya zadacha teorii predstavlenie0 podrazumevaet poisk
ne\-ko\-to\-ro\-go predstavitelya iz fiksirovannogo klassa abstraktnyh
al\-geb\-ra\-i\-ches\-kih struktur (naprimer, algebr Li ili
assotsiativnyh algebr s kvadratichnymi sootnosheniyami) po ego
predstavleniyu, pri e1tom, chto takoe predstavlenie algebraicheskogo
obp2ekta iz zadannogo klassa uzhe opredeleno. Pri reshenii podobnoe0
obratnoe0 zadachi voznikayut, kak pra\-vi\-lo, dostatochno interesnye
novye matematicheskie obp2ekty. Na e1tom puti byli vyyavleny, naprimer,
kvantovye gruppy [4]; nekotorye drugie illyustratsii (algebry
Raka-Vignera, $\mho$--algebry, algebry Sklyanina i t.d.) soderzhat\-sya
v [1,3].

Odnako, obratnaya zadacha teorii predstavlenie0 mozhet bytp1
ras\-smot\-re\-na i v drugom rakurse. Budem predpolagatp1, chto
abstraktnye0 obp2ekt prinadlezhit dostatochno prostomu klassu
struktur, naprimer, yavlyaet\-sya algebroe0 Li, no sam sposob
predstavleniya otlichen ot standartnogo. Podobnaya postanovka
zadachi byla sformulirovana v [5]. Konechno zhe, metodologicheski
perehod k novomu rakursu oznachaet smenu poiska po\-do\-biya mezhdu
konkretnym i abstraktnym obp2ektami (gomomorfnosti pervogo
vtoromu) na analiz vozmozhnosti predstavleniya vtorogo pervym.
E1to predpolagaet tak ili inache vopros o tom, chto takoe
predstavlenie i chto znachit predstavlyatp1. Kak sledstvie,
privlechenie abstraktnogo al\-geb\-ra\-i\-ches\-ko\-go apparata teorii
kategorie0 predstavlyaet\-sya razumnym, chto i bylo sdelano v
[6]. Takim obrazom, dannaya statp1ya, posvyashchennaya tomu, kak
algebra Virasoro (bolee tochno, algebra Vitta, tsentralp1nym
ras\-shi\-re\-ni\-em kotoroe0 algebra Virasoro yavlyaet\-sya) voznikaet
pri primenenii ka\-te\-gor\-noe0 teorii predstavlenie0 k analizu
$q_R$--konformnyh simmetrie0 v dvumernoe0 nelokalp1noe0
$q_R$--konformnoe0 teorii polya, predstavlyaet soboe0
sistematicheskoe izlozhenie materiala e1lektronnyh preprintov
[5,6].

\head \cyb 1. $q_R$--konformnye teorii polya i $q_R$--konformnye
simmetrii [7]
\endhead

Spektr polee0 prostee0shee0 modeli konformnoe0 teorii polya
porozhdaet\-sya odnim svobodnym bozonnym polem i, takim obrazom,
sostoit iz so\-stav\-nyh polee0, poluchennyh pri pomoshchi normalp1no
uporyadochennogo po\-to\-chech\-no\-go proizvedeniya postoennogo po
svobodnomu polyu toka (sm. obshchuyu kon\-st\-ruk\-tsiyu postroeniya toka
po polyu Fubini-Venetsiano v [7]) na sebya ko\-nech\-noe chislo raz.
Ostalp1nye modeli konformnoe0 teorii polya, polnoe opi\-sa\-nie kotoryh
dano v [8], poluchayut\-sya pri perehode k neskolp1kim svo\-bod\-nym
polyam (t.n. predstavlenie svobodnyh polee0 [9]).

Modelp1 s odnim svobodnym bozonnym polem opisyvaet\-sya
lagranzhianom
$$L=\int\varphi(z)\bar\varphi(\bar z)\,dzd\bar z,$$
gde integrirovanie vedet\-sya po vnutrennosti kompleksnogo diska.
S mno\-gih tochek zreniya (kak teoreticheskih, tak i tehnicheskih)
polezno ras\-smat\-ri\-vatp1 nelokalp1nuyu deformatsiyu (sr.[10]) e1toe0
modeli s nelokalp1nym lag\-ran\-zhi\-a\-nom
$$L_h=\int K_h(z,\bar z)\varphi(z)\bar\varphi(\bar z)\,dzd\bar z,$$
gde yadro $K_h(z,\bar z)$ predstavlyaet soboe0 funktsiyu Bergmana
metriki Lo\-ba\-chev\-s\-ko\-go v edinichnom kompleksnom diske:
$$K_h(z,\bar z)=(1-z\bar z)^{-2-2h},$$
v silu chego nelokalp1nye0 lagranzhian okazyvaet\-sya invariantnym
ot\-no\-si\-telp1\-no vseh mebiusovyh (drobno-linee0nyh) preobrazovanie0
kompleksnogo diska. Otmetim, chto lagranzhian svobodnogo bozonnogo
polya $L(\varphi)$ po\-lu\-cha\-et\-sya iz nelokalp1nyh lagranzhianov
$L_h$ predelp1nym perehodom $h\to\infty$. Modeli dvumernoe0
kvantovoe0 teorii polya s lagranzhianom $L_h$ (i es\-tes\-t\-ven\-ny\-mi
nelokalp1nymi skobkami Puassona, ke1lerovym potentsialom otvechayushchee0
kotorym ke1lerovoe0 metriki sluzhit lagranzhian $L_h$) yavlyayut\-sya
prostee0shimi modelyami nelokalp1noe0 teorii, i v to zhe samoe vremya
oni obladayut interesnoe0 i netrivialp1noe0 vnutrennee0 algebraicheskoe0
strukturoe0.

Na yazyke operatornogo formalizma kvantovoe0 teorii polya
ukazannye ne\-lo\-kalp1\-nye modeli posle golomorfno-antigolomorfnoe0
(kiralp1noe0) fak\-to\-ri\-za\-tsii byli rassmotreny v rabote avtora [7]
i nazvany $q_R$--kon\-form\-ny\-mi teoriyami polya ($q_R=\frac1{2h-1}$).
Komponenty operatornyh polee0 $q_R$--konformnyh teorie0 polya, kotorye
invariantny otnositelp1no vseh kom\-p\-lek\-s\-nyh proektivnyh
preobrazovanie0 sfery Rimana, dopuskayut re\-a\-li\-za\-tsiyu tenzornymi
operatorami v modulyah Verma nad algebroe0 Li $\sltwo$ i, kak
sledstvie, zapisyvayut\-sya (obobshchennymi) differentsialp1nymi
ope\-ra\-to\-ra\-mi, poe1tomu mnogie zadachi teorii polya dlya e1tih
modelee0 dopuskayut ana\-li\-ti\-ches\-koe ili chislennoe resheniya, pri
e1tom formulirovka ih raz\-nost\-nyh (reshetochnyh) approksimatsie0 ne
predstavlyaet truda (naprimer, za\-da\-cha integrirovaniya
kvantovo-polevogo volchka, rassmatrivavshayasya dlya konformnyh
teorie0 polya v [11] i dlya $q_R$--konformnyh teorie0 v [12]).

Vazhnoe0 osobennostp1yu dvumernyh konformnyh teorie0 polya
yavlyaet\-sya nalichie beskonechnyh dopolnitelp1nyh simmetrie0,
opisyvaemyh algebroe0 Vitta ili ee tsentralp1nym rasshireniem --
algebroe0 Virasoro $\Cvir$ [13]. E1ti simmetrii poyavlyayut\-sya
kak komponenty razlozheniya (kiralp1no fak\-to\-ri\-zo\-van\-no\-go)
tenzora e1nergii-impulp1sa, kvadrata svobodnogo bozonnogo polya (bolee
tochno, postroennogo po nemu toka). Analogi predstavleniya svobodnyh
polee0 predstvaleniya svobodnyh polee0 i tenzora e1nergii-im\-pulp1\-sa
v $q_R$--konformnom sluchae byli podrobno izucheny v rabote [7].
Tak $q_R$--konformnye0 tenzor e1nergii-impulp1sa yavlyaet\-sya
proizvodyashchee0 funk\-tsi\-ee0 tenzornyh operatorov spina 2 v module
Verma nad algebroe0 Li $\sltwo$. Esli modulp1 Verma $V_h$ ($h$ --
e1kstremalp1nye0 ves [14]) nad e1toe0 algebroe0 realizovan v
prostranstve mnogochlenov $\Bbb C[z]$ ot odnogo kompleksnogo
pe\-re\-men\-no\-go $z$, a dee0stvie generatorov $L_i$ ($i=-1,0,1$;
$[L_i,L_j]=(i-j)L_{i+j}$) algebry Li $\sltwo$ v nem zapizyvaet\-sya
v vide
$$L_{-1}=z,\quad L_0=z\partial_z+h,\quad L_1=z\partial_z^2+2h\partial_z,$$
to komponenty $q_R$--konformnogo tenzora e1nergii-impulp1sa
($q_R$--kon\-form\-nye simmetrii) imeyut vid:
$$ L_{-k}=(\xi+h(k+1))\partial_z^k\quad (k\ge0),\quad
L_k=z^k\frac{\xi+h(k+1)}{(\xi+2h)\ldots(\xi+2h+k-1)}\quad (k\ge1),$$
где $\xi=z\partial_z.$

$q_R$--konformnye simmetrii obrazuyut beskonechnoe semee0stvo. Byl
pred\-pri\-nyat ryad popytok reshitp1 obratnuyu zadachu teorii
predstavlenie0 dlya nih, t.e. vyyavitp1 ih algebraicheskuyu
strukturu. Odna iz popytok opi\-sa\-na v rabote [7], ee rezulp1tatom
yavlyaet\-sya assotsiativnaya algebra $\Cvir(q_R)$, analog
nelinee0noe0 $\sla_2$ raboty [15]. Drugoe0 podhod byl namechen v
e1lek\-t\-ron\-nom preprinte [5] i razvit v sdannoe0 v pechatp1 rabote
avtora ``Pri\-bli\-zhen\-nye predstavleniya i algebra Virasoro''. Kak
by to ni bylo, ostaet\-sya neobhodimostp1 dalp1nee0shego poiska
inyh putee0 resheniya obratnoe0 zadachi teorii predstavlenie0 dlya
$q_R$--konformnyh simmetrie0. Sleduyushchie0 pa\-rag\-raf posvyashchen
opisaniyu odnogo iz vozmozhnyh variantov.

\head \cyb 2. Algebraicheskaya struktura $q_R$--konformnyh simmetrie0
\endhead

\definition{\cyb Opredelenie 1 [5]}

{\cyb A.} \cyr Linee0noe prostranstvo $\Fv$ nazyvaet\-sya {\cyi
lievskim kompozitom}, esli v nem fiksirovany podprostranstva
$\Fv_1,\ldots \Fv_n$ ($\dim\Fv_i>1$), snabzhennye soglasovannymi
strukturami algebr Li. Soglasovannostp1 oznachaet, chto v kazhdom
podprostranstve $\Fv_i\cap\Fv_j$ korrektno opredelena odna i ta zhe
struktura algebry Li, indutsirovannaya iz $\Fv_i$ ili iz $\Fv_j$,
t.e. vlozheniya $\Fv_i\cap\Fv_j$ v $\Fv_i$ i $\Fv_j$ sutp1 vlozheniya
algebr Li. Lievskie0 kompozit nazyvaet\-sya {\cyi plotnym}, esli
$\Fv_1\uplus\ldots\uplus\Fv_n=\Fv$ (zdesp1 $\uplus$ oboznachaet
summu linee0nyh prostranstv). Lievskie0 kompozit nazyvaet\-sya
{\cyi svyaznym}, esli dlya lyubyh $i$ i $j$ sushchestvuet
posledovatelp1nostp1 $k_1,\ldots k_m$ ($k_1=i$, $k_m=j$) takaya, chto
$\Fv_{k_l}\cap\Fv_{k_{l+1}}\ne\varnothing$.

{\cyb B.} \cyr {\cyi Predstavleniem\/} lievskogo kompozita $\Fv$ v
prostranstve $H$ nazyvaet\-sya linee0noe otobrazhenie
$T:\Fv\mapsto\End(H)$ takoe, chto $\left.T\right|_{\Fv_i}$ --
predstavlenie algebry Li $\Fv_i$ pri lyubom vybore $i$.

{\cyb V.} \cyr Pustp1 $\Fg$ -- algebra Li. Linee0noe otobrazhenie
$T:\Fg\mapsto\End(H)$ na\-zy\-va\-et\-sya {\cyi kompozitnym predstavleniem\/}
$\Fg$ v linee0nom prostranstve $H$, es\-li sushchestvuet mnozhestvo
$\Fg_1,\ldots,\Fg_n$ podalgebr Li algebry $\Fg$, obrazuyushchih
plotnye0 svyaznye0 lievskie0 kompozit, i $T$ -- ego predstavlenie.
\enddefinition

\cyr
Privodimostp1 i neprivodimostp1, razlozhimostp1 i nerazlozhimostp1
pred\-stav\-le\-nie0 lievskogo kompozita, a takzhe estestvennye operatsii
nad nimi, op\-re\-de\-lya\-yut\-sya tak zhe kak i dlya algebr Li [16].
Mozhno takzhe sfor\-mu\-li\-ro\-vatp1 superanalog opredeleniya 1.
Sovokupnostp1 predstavlenie0 fik\-si\-ro\-van\-no\-go lievskogo kompozita
zamknuto otnositelp1no tenzornogo pro\-iz\-ve\-de\-niya i, sledovatelp1no,
dopuskaet strukturu {\cyi tenzornoe0 kategorii} [16].

\example{\cyi Primer 1 (Oktae1dricheskie0 lievskie0 kompozit [5])}
\cyr Rassmotrim oktae1dr s vershinami $A$, $B$, $C$, $D$, $E$, $F$,
rebrami $(AB)$, $(AC)$, $(AD)$, $(AE)$, $(BC)$, $(BF)$, $(CD)$,
$(CF)$, $(DE)$, $(DF)$, $(EF)$, i granyami $(ABC)$, $(ACD)$, $(ADE)$,
$(AEB)$, $(BCF)$, $(CDF)$, $(DEF)$, $(EBF)$. Pustp1 $\Fv$ --
shestimernoe linee0noe prostranstvo s bazisom, zanumerovannym
vershinami oktae1dra, $\Fv_1$, $\Fv_2$, $\Fv_3$, $\Fv_4$ --
chetverka trehmernyh podprostranstv v $\Fv$, otvechayushchih granyam
$(ABC)$, $(ADE)$, $(CDF)$, $(EBF)$. Vse podprostranstva $\Fv_i$ snabdim
strukturami algebr Li, izomorfnyh algebre $\sothree$ (e1ti struktury
dolzhny bytp1 so\-gla\-so\-va\-ny s orientatsiyami granee0). Pyaterka
$(\Fv,\Fv_1,\Fv_2,\Fv_3,\Fv_4)$ yavlyaet\-sya plot\-nym svyaznym
kompozitom.

\proclaim{\cyb Predlozhenie} \cyi Pustp1 $T$ -- proizvolp1noe
predstavlenie lievskogo kom\-po\-zi\-ta $(\Fv;\Fv_1,\Fv_2,\Fv_3,\Fv_4)$ v
konechnomernom linee0nom prostranstve $H$, togda $H$ nadelyaet\-sya
strukturoe0 predstavleniya algebry Li $\sofour$. Esli $T$ --
ne\-pri\-vo\-di\-moe predstavlenie, to sushchestvuyut veshchestvennye chisla
$\lambda_A$, $\lambda_B$, $\lambda_C$, $\lambda_D$, $\lambda_E$, $\lambda_F$
takie, chto operatory $T(A)-\lambda_A\boldkey 1$,
$T(B)-\lambda_B\boldkey 1$, $T(C)-\lambda_C\boldkey 1$,
$T(D)-\lambda_D\boldkey 1$, $T(E)-\lambda_E\boldkey 1$,
$T(F)-\lambda_F\boldkey 1$ obrazuyut neprivodimoe predstavlenie $\sofour$.
\endproclaim

\demo{\cyi Dokazatelp1stvo} \cyr Vo-pervyh, otmetim, chto kommutator
operatorov, ot\-ve\-cha\-yu\-shchih protivopolozhnym vershinam, kommutiruet
s operatorami, ot\-ve\-cha\-yu\-shchi\-mi chetyrem drugim vershinam. Kak
sledstvie, on kommutiruet so vsemi shestp1yu operatorami,
predstavlyayushchimi generatory ok\-ta\-e1d\-ri\-ches\-ko\-go lievskogo
kompozita, poskolp1ku oni vyrazhayut\-sya cherez kommutatory
upomyanytyh vyshe chetyreh operatorov. Itak, kommutatory
operatorov, otvechayushchih protivopolozhnym vershinam oktae1dra,
lezhat v tsentre al\-geb\-ry Li $\Fg$, porozhdennoe0 vsemi shestp1yu
operatorami, predstavlyayushchimi generatory kompozita.
Otfaktorizuem algebru Li $\Fg$ po tsentru. Po\-lu\-chen\-naya faktoralgebra
izomorfna $\sofour$ (pri e1tom ispolp1zuet\-sya to, chto for\-mu\-ly
dlya kommutatorov vseh shesti operatorov izvestny s tochnostp1yu
do e1lementov tsentra algebry Li $\Fg$). Utverzhdenie predlozheniya
estp1 sled\-s\-t\-vie e1togo rezulp1tata i togo fakta, chto lyuboe
tsentralp1noe ras\-shi\-re\-nie poluprostoe0 algebry Li trivialp1no
(t.e. rasshchepimo -- sm.napr.[17])\qed
\enddemo
\endexample

\cyr
Konstruktsiya oktae1dricheskogo lievskogo kompozita obobshchaet\-sya
na ne\-ko\-to\-rye0 estestvennye0 klass polie1drov, odnako, v obshchem
sluchae analog do\-ka\-zan\-no\-go vyshe predlozheniya neizvesten.

\example{\cyi Primer 2 (Lievskie0 kompozit Vitta [5])} \cyr Pustp1
$\Fw$ -- algebra Vitta [18] (podalgebra kompleksifikatsii $\CVect$
algebry Li $\Vect$ gladkih vektornyh polee0 na okruzhnosti $\crc$,
porozhdennaya vsemi loranovskimi polinomialp1nymi vektornymi polyami)
s bazisom $e_k$ ($k\in\Bbb Z$) i kom\-mu\-ta\-tsi\-on\-ny\-mi sootnosheniyami
$[e_i,e_j]=(i-j)e_{i+j}$. Algebra Virasoro $\Cvir$ [19] yavlyaet\-sya
odnomernym netrivialp1nym tsentralp1nym rasshireniem algebry Vitta.

Rassmotrim dve podalgebry $\Fp_{\pm}$ algebry Li $\Fw$, porozhdennye
$e_i$ s $i\ge-1$ i $i\le1$; otmetim, chto $\Fp_+\cap\Fp_-=\sltwo$.
Troe0ka $(\Fw;\Fp_+,\Fp_-)$ yavlyaet\-sya plotnym svyaznym lievskim
kompozitom, kotorye0 my budem nazyvatp1 kompozitom Vitta.

Lyuboe predstavlenie algebry Virasoro zadaet predstavlenie li\-ev\-s\-ko\-go
kompozita Vitta v tom zhe prostranstve. Obratnoe, konechno zhe,
neverno, chto pokazyvaet sleduyushchaya teorema.

\proclaim{\cyb Teorema 1A} \cyi Dee0stvie algebry Li $\sltwo$ v
lyubom module Verma $V_h$ ($h$ -- e1kstremalp1nye0 ves) mozhet bytp1
odnoznachno prodolzheno do predstavleniya lievskogo kompozita Vitta
$(\Fw;\Fp_+,\Fp_-)$ i, sledovatelp1no, do kompozitnogo predstavleniya
algebry Vitta $\Fw$ (i algebry Virasoro $\Cvir$).
\endproclaim

\demo{\cyi Dokazatelp1stvo} \cyr Utverzhdenie teoremy sleduet iz
yavnyh formul dlya ten\-zor\-nyh operatorov v module Verma nad algebroe0
Li $\sltwo$ spina 2 ($q_R$--konformnyh simmetrie0), kotorye i
osushchestvlyayut predstavlenie li\-ev\-s\-ko\-go kompozita Vitta.
\enddemo
\endexample

\remark{\cyi Zamechanie 1} \cyr Konstruktsiya lievskogo kompozita Vitta
obobshchaet\-sya na ri\-ma\-no\-vy poverhnosti starshih rodov v duhe
I.M.Krichevera i S.P.No\-vi\-ko\-va [19].
\endremark

\remark{\cyi Zamechanie 2} \cyr Obobshchaya terminologiyu rabot [1,3]
mozhno govoritp1, chto tenzornye operatory spina 2 (t.e.
$q_R$--konformnye simmetrii) v modulyah Verma $V_h$ nad algebroe0 Li
$\sltwo$ obrazuyut semee0stvo skrytyh sim\-met\-rie0, chp1ya
algebraicheskaya struktura estp1 struktura lievskogo kompozita Vitta.
\endremark

\cyr Otmetim, chto esli skrytye simmetrii realizuyut predstavlenie
li\-ev\-s\-ko\-go kompozita, oni ne trebuyut ``raspakovyvaniya'' (podobnaya
situatsiya imeet mesto v sluchae izokommutatornyh algebr skrytyh
simmetrie0 i lievskih $\Fg$--puchkov [3:Topic 3;1:\S2.1]).

\example{\cyi Primer 3} \cyr Pustp1 $\Fw$ -- algebra Vitta i
$(\Fw;\Fp_{\pm})$ -- lievskie0 kompozit Vitta. Rassmotrim abelevo
(netsentralp1noe) rasshirenie $\Fw^e$ algebry Vitta pri pomoshchi
generatorov $f_i$ ($i\in\Bbb Z$) takoe, chto $[e_i,f_j]=jf_j$.
Podalgebry $\Fp_{\pm}$ algebry Li $\Fw$ mogut bytp1 rasshireny do
podalgebr $\Fp^e_{\pm}$ algebry Li $\Fw^e$ generatorami $f_i$,
gde $i\ge 0$ i $i\le 0$, sootvet\-stvenno. Troe0ka $(\Fw^e;\Fp^e_{\pm})$
zadaet rasshirennye0 kompozit Vitta.

\proclaim{\cyb Teorema 1B} \cyi Predstavlenie kompozita Vitta v
lyubom module Verma $V_h$ ($h$ -- e1kstremalp1nye0 ves) nad
algebroe0 Li $\sltwo$ mozhet bytp1 prodolzheno do predstavleniya
rasshirennogo kompozita Vitta $(\Fw^e;\Fp^e_{\pm})$ i,
sle\-do\-va\-telp1\-no, kompozitnoe predstavlenie algebry Vitta $\Fw$ v
$V_h$ prodolzhaet\-sya do kompozitnogo predstavleniya $\Fw^e$.
\endproclaim

\demo{\cyi Dokazatelp1stvo} \cyr Dopolnitelp1nye generatory $f_i$
predstavlyayut\-sya ten\-zor\-ny\-mi operatorami spina 1 (komponentami
razlozheniya $q_R$--affinnogo toka [7]), a imenno
$f_i\mapsto\partial_z^i$ ($i\ge 0$),
$f_{-i}\mapsto z^i\frac1{(\xi+2h)\ldots(\xi+2h+i-1)}$ ($i\ge 1$)
\qed\enddemo
\endexample

\remark{\cyi Zamechanie 3} \cyr Vesp1ma interesno bylo by
rassmotretp1 kompozitnye preds\-tav\-le\-niya veshchestvennyh poluprostyh
algebr Li, ``sostavlennye'' iz obych\-nyh predstavlenie0 nekotoryh
estestvennyh podalgebr (naprimer, dvuh pro\-ti\-vo\-po\-lozh\-nyh
maksimalp1nyh parabolicheskih podalgebr ili dvuh
pro\-ti\-vo\-po\-lozh\-nyh borelevskih algebr plyus, vozmozhno, nekotoryh
vlozhennyh algebr $\sltwo$, i t.d.).
\endremark

\remark{\cyi Zamechanie 4} \cyr Kompozitnye predstavleniya nekotoryh
drugih al\-geb\-ra\-i\-ches\-kih struktur (izotopicheskih par)
rassmatrivalisp1 v rabote avtora [20].
\endremark

Itak, vo vtorom paragrafe my vyyasnili, chto $q_R$--konformnye
sim\-met\-rii mozhno schitatp1 realizatsiee0 toe0 zhe samoe0
algebraicheskoe0 struktury (algebry Vitta ili algebry Virasoro),
chto i konformnye simmetrii, od\-na\-ko, pri nekotorom novom ponimanii
samoe0 operatsii ``pred\-stav\-le\-niya''. Takim obrazom, estestvenno
postavitp1 vopros o tom, chto znachit ``predstavlyatp1'' i kakimi
svoe0stvami dolzhny obladatp1 nestandartnye predstavleniya, chtoby s
nimi mozhno bylo rabotatp1 takzhe kak s obychnymi, a takzhe kak
poluchatp1 iz uzhe imeyushchihsya tipov predstavlenie0 novye
ne\-iz\-ves\-t\-nye tipy. Na e1ti voprosy prizvan otvetitp1 sleduyushchie0
paragraf, kotorye0 budet nositp1 dostatochno formalp1nye0
matematicheskie0 harakter, dalekie0 ot nachalp1noe0 postanovki
zadachi iz nelokalp1noe0 kvantovoe0 te\-o\-rii polya.

\head\cyb 3. Kategornaya teoriya predstavlenie0 i $q_R$--konformnye
simmetrii: kompozitnye i overlee0nye predstavleniya
\endhead

\subhead\cyb 3.1. E1lementy kategornoe0 teorii predstavlenie0 [6]
\endsubhead \cyr Dalee my budem rassmatrivatp1 predstavleniya
klassov obp2ektov, obrazuyushchih kategoriyu, kotoruyu my
budem nazyvatp1 bazovoe0 kategoriee0. Kategornye aspekty
standartnoe0 teorii predstavlenie0 razbiralisp1 v [16].
Nekotorye ka\-te\-gor\-nye obobshcheniya opisany v [21]. Odnako, nizhe
budut sformulirovan naibolee obshchie0 abstraktnye0 formalizm,
kotorye0 neobhodim dlya nashih zadach.

\definition{\cyb Opredelenie 2A} \cyr {\cyi Teoriee0 predstavlenie0\/}
bazovoe0 kategorii $\CA$ na\-zy\-va\-et\-sya kontravariantnye0 funktor $R$
iz kategorii $\CA$ v kategoriyu $\Abel$ vseh malyh abelevyh
kategorie0.
\enddefinition

\cyr Inogda sleduet rassmatrivatp1 kategoriyu $\Add$ vseh malyh
additivnyh kategorie0 vmesto $\Abel$. Odnako, dlya prostoty budet
ispolp1zovatp1sya tolp1ko kategoriya $\Abel$.

Chasto bazovaya kategoriya imeet ryad horoshih svoe0stv, naprimer,
chto dlya proizvolp1nogo konechnogo semee0stva obp2ektov
sushchestvuet ih ko\-pro\-iz\-ve\-de\-nie, sovpadayushchee s proizvedeniem.
E1ta situatsiya imeet mesto dlya algebr Li, grupp Li, konechnyh
grupp, assotsiativnyh algebr, algebr Hopfa i mnogih drugih struktur.
Odnako, bolp1shinstvo algebraicheskih par (naprimer, izotopicheskie
pary [1:\S2.2;19]) ne obrazuyut kategorie0 podobnogo tipa. Ot
bazovoe0 kategorii $\CA$, v kotoroe0 proizvedeniya i
ko\-pro\-iz\-ve\-de\-niya konechnogo chisla obp2ektov sushchestvuyut
i sovpadayut, my po\-tre\-bu\-em v opredelenii teorii predstavlenie0,
chtoby bylo zadano as\-so\-tsi\-a\-tiv\-noe semee0stvo vlozhenie0
$R(a)\times R(b)\hookrightarrow R(a+b)$ ($a$ i $b$ --
proizvolp1nye obp2ekty bazovoe0 kategorii $\CA$). Takie teorii
predstavlenie0 budem na\-zy\-vatp1 {\cyi kvazitenzorialp1nymi}.

\remark{\cyi Zamechanie 5} \cyr Esli obp2ekt $a$ bazovoe0 kategorii
$\CA$ dopuskaet ko\-as\-so\-tsi\-a\-tiv\-nye0 monomorfizm $\varepsilon$
v $a+a$, to $R(a)$ yavlyaet\-sya tenzornoe0 kategoriee0, kolp1 skoro
teoriya predstavlenie0 $R$ kvazitenzorialp1na.
\endremark

Sushchestvuyut nekvazitenzorialp1nye teorii predstavlenie0 dazhe
dlya ho\-ro\-sho izvestnyh kategorie0 predstavlyaemyh obp2ektov,
naprimer, obshchie $\HS$--proektivnye predstavleniya algebr Li [5]
ili unitarnye $\HS$--psev\-do\-pred\-s\-tav\-le\-niya grupp Li.

\definition{\cyb Opredelenie 2B} \cyr Teoriya predstavlenie0 $R$ dlya
bazovoe0 kategorii $\CA$ na\-zy\-va\-et\-sya {\cyi gomomorfnoe0}, esli
sushchestvuet podkategoriya $\CA_0$ kategorii $\CA$ ({\cyi podkategoriya
mishenee0\/}) takaya, chto dlya proizvolp1nogo obp2ekta $a$
ka\-te\-go\-rii $\CA$ kategoriya $R(a)$ otozhdestvlyaet\-sya
s kategoriee0 $\Mor(a,\CA_0)$ vseh (klassov e1kvivalentnosti)
morfizmov iz $a$ v obp2ekty kategorii mishenee0 $\CA_0$.
\enddefinition

\cyr Naprimer, teorii vseh linee0nyh, proektivnyh ili unitarnyh
pred\-s\-tav\-le\-nie0 grupp Li gomomorfny. Otmetim, chto kategoriya
mishenee0 $\CA_0$ vsegda yavlyaet\-sya additivnoe0 (no ne
obyazatelp1no abelevoe0) podkategoriee0 bazovoe0 kategorii $\CA$.

\definition{\cyb Opredelenie 2V} \cyr Teoriya predstavlenie0 $R$
dlya bazovoe0 kategorii $\CA$ na\-zy\-va\-et\-sya {\cyi skryto
gomomorfnoe0}, esli sushchestvuet gomomorfnaya teoriya pred\-s\-tav\-le\-nie0
$R'$ dlya kategorii $\CK$ i funktor (kak pravilo, mnogoznachnye0)
$\varrho:\CA\mapsto\CK$ takie, chto
$R=R'\circ\varrho$.
\enddefinition

\cyr Nizhe my privedem nekotorye primery i obshchie konstruktsii
skryto gomomorfnyh teorie0 predstavlenie0 dlya kategorii $\Lie$
algebr Li, in\-s\-pi\-ri\-ro\-van\-nye predlozhennym vyshe resheniem
obratnoe0 zadachi teorii predstavlenie0 dlya $q_R$--konformnyh simmetrie0 v
nelokalp1noe0 kvantovoe0 teorii polya, kotorye ne yavlyayut\-sya
gomomorfnymi, i opishem ih in\-ter\-p\-re\-ta\-tsii v terminah kategornoe0
teorii predstavlenie0.

\subhead\cyb 3.2. Kompozitnye teorii predstavlenie0 [6]
\endsubhead \cyr Izlozhim abstraktnye0 kategornye0 formalizm
dlya kompozitnyh predstavlenie0 algebr Li, op\-re\-de\-len\-nyh vo vtorom
paragrafe.

\definition{\cyb Opredelenie 3A} \cyr Pustp1 $\CA$ --
topologizirovannaya bazovaya kategoriya (t.e. snabzhennaya topologiee0
Grotendika [22]). Pustp1 $R$ -- teoriya pred\-s\-tav\-le\-nie0 dlya $\CA$.
{\cyi Kompozitnaya teoriya predstavlenie0\/} $\BC(R)$ dlya $\CA$ mozhet
bytp1 postroena sleduyushchim obrazom. Pustp1 $a$ -- obp2ekt
bazovoe0 kategorii $\CA$ i $S=(s_1,s_2,\ldots s_n)$ ($s_i\in\Mor(a_i,a)$)
 -- pokrytie $a$, togda obp2ekty ka\-te\-go\-rii $\BC(R)(a)$ sutp1 nabory
dannyh $(b_1,b_2,\ldots b_n)$, $b_i\in R(a_i)$ takih, chto dlya
proizvolp1nogo obp2ekta $c$ i monomorfizmov $f\in\Mor(c,a)$,
$f_i\in\Mor(c,a_i)$ ($f=s_i\circ f_i$) vypolnyaet\-sya ravenstvo
({\cyi kompozitnoe pravilo sklee0ki\/})
$$R(f_i)^*(b_i)=R(f_j)^*(b_j).$$
Morfizmy v $\BC(R)(a)$ opredelyayut\-sya soglasovannym obrazom.
\enddefinition

\cyr
Dlya lyuboe0 teorii predstavlenie0 $R$ kompozitnaya teoriya predstavlenie0
$\BC(R)$ yavlyaet\-sya puchkom abelevyh kategorie0 nad
topologizirovannoe0 ba\-zo\-voe0 kategoriee0 $\CA$ [22]. E1to -- puchok,
kanonicheski postroennye0 po pred\-puch\-ku $R$ nad topologizirovannoe0
bazovoe0 kategoriee0 $\CA$ (otmetim, chto teoriya predstavlenie0 dlya
topologizirovannoe0 bazovoe0 kategorii $\CA$ estp1 v tochnosti po
opredeleniyu predpuchok nad nee0).

\proclaim{\cyb Teorema 2A [6]} \cyi Kompozitnye predstavleniya
algebr Li obrazuyut kom\-po\-zit\-nuyu teoriyu predstavlenie0 $\BC(R)$,
gde $R$ -- standartnaya teoriya pred\-s\-tav\-le\-nie0 algebr Li (nakrytiya
algebr Li zadayut\-sya kak plotnye svyaznye lievskie kompozity).
\endproclaim

\cyr Otmetim, chto topologiya Grotendika v teoreme 2A otlichaet\-sya
ot obych\-noe0 [21].

\remark{\cyi Zamechanie 6} \cyr Esli $R$ -- standartnaya teoriya
predstavlenie0 algebr Li, to kompozitnaya teoriya predstavlenie0
$\BC(R)$ skryto gomomorfna, vspo\-mo\-ga\-telp1\-naya kategoriya $\CK$ --
kategoriya lievskih kompozitov, a podkategoriya mishenee0 $\CK_0$
sostoit iz algebr Li $\End(H)$ linee0nyh operatorov v li\-nee0\-nyh
prostranstvah, t.e. ta zhe, chto i dlya standartnoe0 gomomorfnoe0
teorii predstavlenie0. Odnako, esli $R$ -- obshchaya teoriya
predstavlenie0, to o skrytoe0 gomomorfnosti $\BC(R)$ nichego ne izvestno.
\endremark

\cyr Ya podozrevayu, chto ponyatie skrytoe0 gomomorfnosti
kompozitnyh te\-o\-rie0 predstavlenie0 mozhet bytp1 ponyato v
terminah teorii toposov [22].

\remark{\cyi Zamechanie 7} $\BC(\BC(R))=\BC(R)$.
\endremark

\subhead\cyb 3.3. Overlee0nye teorii predstavlenie0 [6]
\endsubhead \cyr Ogranichennostp1 kom\-po\-zit\-noe0 teorii predstavlenie0
polnostp1yu vyyavlyaet\-sya na primere kom\-po\-zit\-nyh predstavlenie0
algebry Vitta skrytymi beskonechnomernymi $q_R$--kon\-form\-ny\-mi
simmetriyami v modulyah Verma nad algebroe0 Li $\sltwo$.
Vo-per\-vyh, tenzornoe proizvedenie lyubogo chisla e1tih
neprivodimyh kom\-po\-zit\-nyh predstavlenie0 snova neprivodimo.
E1tot fakt rezko protivorechit intuitsii. Vo-vtoryh, skrytye
simmetrii ne obrazuyut sami nikakogo kompozitnogo predstavleniya
algebry Vitta, v to vremya kak intuitivno oni dolzhny
obrazovyvatp1 ``prisoedinennoe'' predstavlenie. S tochki zre\-niya
matematicheskoe0 fiziki perechislennye nedostatki sushchestvenny,
po\-skolp1\-ku oni ne pozvolyayut ispolp1zovatp1 v polnoe0 mere
apparat teorii predstavlenie0 dlya analiza protsessov raspada,
sinteza i vzaimodee0stviya kvazichastits v dvumernoe0 nelokalp1noe0
teorii polya.

Takim obrazom, vyyavlyaet\-sya neobhodimostp1 obobshcheniya ponyatiya
kom\-po\-zit\-no\-go predstavleniya. S e1toe0 tselp1yu opredelim operatornye
lievskie kompozity $\LC(H)$ kak nabory prostranstv operatorov
$\End(H_i)$, pod\-pros\-t\-ranstv prostranstva operatorov $\End(H)$
($H=H_1+\ldots+H_m$), s es\-tes\-t\-ven\-ny\-mi strukturami algebr Li.

\definition{\cyb Opredelenie 4}

{\cyb A.} \cyr {\cyi Overlee0nym predstavleniem\/} lievskogo
kompozita $\Fv$ v prostranstve $H$ nazyvaet\-sya gomomorfizm
$T$ kompozita $\Fv$ v operatornye0 lievskie0 kom\-po\-zit $\LC(H)$.

{\cyb B.} \cyr Pustp1 $\Fg$ -- algebra Li. Linee0noe otobrazhenie
$T:\Fg\mapsto\End(H)$ nazyvaet\-sya {\cyi overlee0nym predstavleniem\/}
al\-geb\-ry Li $\Fg$ v linee0nom pros\-t\-ran\-s\-t\-ve $H$, esli
sushchestvuet mnozhestvo $\Fg_1,\Fg_2,\ldots,\Fg_n$ podalgebr Li
al\-geb\-ry $\Fg$, obrazuyushchih plotnye0 svyaznye0 kompozit, i $T$ -- ego
overlee0noe predstavlenie.
\enddefinition

\remark{\cyi Zamechanie 8} \cyr Overlee0nye predstavleniya
proizvolp1noe0 algebry Li $\Fg$ ob\-ra\-zu\-yut tenzornuyu kategoriyu.
\endremark

\cyr Overlee0nye predstavleniya reshayut perechislennye ranee
trudnosti kom\-po\-zit\-noe0 teorii predstavlenie0.

\proclaim{\cyb Teorema 1V [6]} \cyi Tenzornye operatory spina 2 v
modulyah Verma $V_h$ nad algebroe0 Li $\sltwo$ obrazuyut overlee0nye
predstavleniya algebry Vitta, yavlyayushchiesya podpredstavleniyami
$\End(V_h)$.
\endproclaim

\remark{\cyi Zamechanie 9} \cyr Tenzornye operatory lyubogo
neotritsatelp1nogo tselogo spi\-na $n$ v module Verma $V_h$ nad
algebroe0 Li $\sltwo$ (opisannye polnostp1yu v [7]) obrazuyut
overlee0nye predstavleniya algebry Vitta, yavlyayushchiesya
podpredstavleniyami $\End(V_h)$.
\endremark

\cyr
Izlozhim estestvennye0 kategornye0 formalizm dlya konstruktsii
over\-lee0\-nyh predstavlenie0 algebr Li.

\definition{\cyb Opredelenie 3B} \cyr Pustp1 $\CA$ --
topologizirovannaya bazovaya kategoriya. Pustp1 $R$ -- gomomorfhaya
teoriya predstavlenie0 dlya $\CA$ s kategoriee0 mi\-she\-nee0 $\CA_0$,
snabzhennoe0 topologiee0 Grotendika, indutsirovannoe0 iz $\CA$.
{\cyi Over\-lee0\-naya teoriya predstavlenie0\/} $\BO(R)$ dlya $\CA$
mozhet bytp1 postroena sle\-du\-yu\-shchim obrazom. Pustp1 $a$ -- obp2ekt
bazovoe0 kategorii $\CA$ i $S=(s_1,s_2,\ldots s_n)$ ($s_i\in\Mor(a_i,a)$)
-- pokrytie $a$, togda obp2ekty kategorii $\BO(R)(a)$ sostoyat iz
naborov dannyh $(r_1,r_2,\ldots r_n)$, $r_i\in\Mor(a,b_i)$, $b_i$ --
obp2ekty podkategorii mishenee0 $\CA_0$, zadayushchih pokrytie
obp2ekta $b$ toe0 zhe podkategorii mo\-no\-mor\-fiz\-ma\-mi
$t_i\in\Mor(b_i,b)$, tak chto dlya lyubogo podobp2ekta $(c;p]$ obp2ekta
$a$ ($p\in\Mor(c,a)$) vypolnyaet\-sya ravenstvo ({\cyi overlee0noe
pravilo sklee0ki\/}):
$$r_i((a_i;s_i]\cap(c;p])\cap(b_j;t_j]=
(b_i;t_i]\cap r_j((a_j;s_j]\cap(c;p]).$$
Morfizmy v $\BO(R)(a)$ opredelyayut\-sya soglasovannym obrazom.
\enddefinition

\cyr
Kak by to ni bulo, mne ne izvestny opredeleniya overlee0noe0
teorii predstavlenie0 $\BO(R)$ dlya teorii predstavlenie0 $R$,
kotoraya ne yavlyaet\-sya gomomorfnoe0. Otmetim, chto $\BO(R)$
ne yavlyaet\-sya puchkom abelevyh kategorie0 nad $\CA$ v obshchem
sluchae, i ya ne znayu abstraktnogo opisaniya overlee0nyh teorie0
predstavlenie0 v terminah teorii puchkov.

\remark{\cyi Zamechanie 10} \cyr Overlee0nye teorii predstavlenie0
$\BO(R)$, opredelennye dlya gomomorfnyh teorie0 predstavlenie0
$R$, yavlyayut\-sya skryto gomomorfnymi.
\endremark

\proclaim{\cyb Teorema 2B [6]} \cyi Overlee0nye predstavleniya
algebr Li obrazuyut over\-lee0\-nuyu teoriyu predstavlenie0 $\BO(R)$,
gde $R$ -- standartnaya teoriya pred\-s\-tav\-le\-nie0 algebr Li (pokrytiya
algebr Li zadayut\-sya kak plotnye svyaznye lievskie kompozity, a
kategoriya mishenee0 $\CA_0$ sostoit iz vseh algebr Li $\End(H)$).
\endproclaim

\cyr
Otmetim, chto topologiya Grotendika v formulirovke Teoremy 2B ot\-lich\-na
ot standartnoe0.

\remark{\cyi Zamechanie 11} \cyr (A) Esli $R$ -- gomomorfnaya teoriya
predstavlenie0 dlya ba\-zo\-voe0 kategorii $\CA$, to dlya lyubogo obp2ekta
$a$ kategorii $\CA$ kategoriya $\BC(R)(a)$ yavlyaet\-sya podkategoriee0
kategorii $\BO(R)(a)$. (B) $\BC(\BO(R))=\BO(R)$.
\endremark

\cyr
Itak, v rabote rassmotreny simmetrie0nye aspekty dvumernoe0
ne\-lo\-kalp1\-noe0 teorii polya, yavlyayushchee0sya prostee0shee0
deformatsiee0 kon\-for\-m\-no-in\-va\-ri\-an\-t\-noe0 teorii s odnim
svobodnym bozonnym polem (\S1). Reshena ob\-rat\-naya zadacha teorii
predstavlenie0 dlya $q_R$--konformnyh simmetrie0,
yav\-lya\-yu\-shchih\-sya beskonechnomernymi skrytymi simmetriyami
ukazannoe0 te\-o\-rii (\S2). Na baze abstraktnoe0 kategornoe0 teorii
predstavlenie0 razvit matematicheskie0 apparat, dostatochnye0 dlya
primeneniya simmetrie0nyh me\-to\-dov k izucheniyu protsessov sinteza,
raspada i vzaimodee0stviya kva\-zi\-chas\-tits v modelyah dvumernoe0
nelokalp1noe0 kvantovoe0 teorii polya (\S3).

\Refs\nofrills{\cyb Spisok literatury}
\roster
\item" [1]" {\cyie Yurp1ev D.V.}, {\cyre E1kskurs v obratnuyu zadachu teorii
predstavlenie0:} mp\_arc/96-477.
\item" [2]" {\cyie Manin Yu.I.}, {\cyre Kalibrovochnye polya i kompleksnaya
geometriya. M., Nauka, 1986}.
\item" [3]" {\it Juriev D.}, Topics in hidden symmetries. I-IV,VI:
hep-th/9405050, q-alg/9610026, q-alg/9611003, q-alg/9611019, q-alg/9708028;
On the infinite-dimensional hidden sym\-met\-ries. I-III: funct-an/9612004,
funct-an/9701009, funct-an/9702002.
\item" [4]" {\cyie Reshetihin N.Yu., Tahtadzhyan L.A., Faddeev L.D.} /\!/
{\cyre Algebra i anal. 1989. T.1, vyp.1. S.178-206}; {\cyie Drinfelp1d V.G.}
/\!/ {\cyre Zapiski nauchn.semin.LOMI, 1986. T.155. S.19-49}; {\cyie Isaev
A.P.} /\!/ {\cyre E1ChAYa. 1995. T.26, vyp.5. S.1204-1263}.
\item" [5]" {\it Juriev D.}, Topics in hidden symmetries. V:
funct-an/9611003.
\item" [6]" {\it Juriev D.}, Hidden symmetries and categorical
representation theory: q-alg/9612026.
\item" [7]" {\cyie Yurp1ev D.V.} /\!/ {\cyre TMF. 1994. T.101, vyp.3.
S.331-348}.
\item" [8]" {\cyie Yurp1ev D.V.} /\!/ {\cyre Algebra i anal. 1991. T.3, vyp.3.
S.197-205.}
\item" [9]" {\cyie Morozov A.Yu.} /\!/ {\cyre E1ChAYa. 1992. T.23, vyp.1.
S.174-238.}
\item"[10]" {\cyie Nelokalp1nye}, {\cyre nelinee0nye i nerenormiruemye
teorii polya. Dubna, OIYaI, 1970.}
\item"[11]" {\it Alekseev A., Faddeev L., Semenov-Tian-Shansky M.} /\!/
Commun.Math.Phys. 1992. V.149. P.335-345.
\item"[12]" {\cyie Yurp1ev D.V.} /\!/ {\cyre TMF. 1996. T.106, vyp.2.
S.333-352.}
\item"[13]" {\it Belavin A.A., Polyakov A.M., Zamolodchikov A.B.} /\!/
Nucl.Phys.B. 1984. P.241. P.333-380.
\item"[14]" {\cyie Zhelobenko D.P.}, {\cyre Predstavleniya reduktivnyh
algebr Li. M., Nauka, 1994.}
\item"[15]" {\it Ro\v cek M.} /\!/ Phys. Lett. B. 1991. V.255. P.554-557.
\item"[16]" {\cyie Kirillov A.A.}, {\cyre E1lementy teorii predstavlenie0.
M., Nauka, 1978; Vvedenie v teoriyu predstavlenie0 i nekommutativnye0
garmonicheskie0 analiz / Sovrem. probl.matematiki. Fundam.napravleniya 22.
M., VINITI, 1988}.
\item"[17]" {\cyie Gisharde A.}, {\cyre Kogomologii topologicheskih grupp
i algebr Li. M., Mir, 1984.}
\item"[18]" {\cyie Fuks D.B.}, {\cyre Kogomologii beskonechnomernyh algebr
Li. M., Nauka, 1984.}
\item"[19]" {\cyie Krichever I.M., Novikov S.P.} /\!/ {\cyre Funkts.anal.i
ego prilozh. 1987. T.21, vyp.2. S.46-63, T.24, vyp.4. S.47-61, 1989. T.23,
vyp.1. S.1-14}; {\it Krichever I.M., Novikov S.P.} /\!/ J.Geom.Phys. 1988.
V.5. S.631-661 [reprinted in ``Geometry and
physics. Essays in honour of I.M.Gelfand'', Eds. S.Gindikin and I.M.Singer,
Pitagora Editrice, Bologna and Elsevier Sci.Publ., Amsterdam, 1991].
\item"[20]" {\cyie Yurp1ev D.V.} /\!/ {\cyre TMF. 1995. T.105, vyp.1. S.18-28;
1997. T.111, vyp.1. S.149-158}.
\item"[21]" {\cyie Geronimus A.Yu.} /\!/ {\cyre Funkts.anal.i ego prilozh.
1971. T.5, vyp.3. S.22-31}.
\item"[22]" {\it Th\'eorie\/} des topos et cohomologie \'etale des sch\'emas.
Berlin, Heidelberg, New York, 1972-1973; {\it Cohomologie\/} \'etale. Berlin,
Heidelberg, New York, 1977.
\endroster
\endRefs
\enddocument